\def\titlep{Automata computation of branching laws for endomorphisms of Cuntz algebras}
\documentclass[11pt]{article}
\usepackage{graphicx,ifthen}
\usepackage{amsmath}
\usepackage{amscd}
\usepackage{amsthm}
\usepackage{amssymb}

\font\germ=eufm10 at12pt

\def\goth#1{\hbox{\germ#1}}

\newcommand{\qedh}{\hfill\qed \\}
\newcommand{\vv}{\vspace{.3in}}

\newtheorem{Thm}{Theorem}[section]

\newtheorem{defi}[Thm]{Definition}
\newtheorem{lem}[Thm]{Lemma}

\newtheorem{prop}[Thm]{Proposition}

\newcommand{\ww}{\vv\noindent}

\def\cal#1{\mathcal #1}
\def\con{{\cal O}_{N}}
\def\edot{=1,\ldots,N}
\def\pr{\noindent{\it Proof.}\,}

\def\nset#1{\{1,\ldots,N\}^{#1}}
\def\co#1{{\cal O}_{#1}}
\def\ltn{l_{2}({\bf N})}
\def\mms{\mathsf{M}_{\sigma}}
\def\bfsnl{{\rm BFS}_{N}(\Lambda)}
\def\enda{{\rm End}{\cal A}}
\def\endcon{{\rm End}\con}

\addtocounter{footnote}{1}

\def\sftt#1{\setcounter{equation}{0}\addtocounter{footnote}{1}\section{#1}}
\def\ssft#1{\subsection{#1}}

\begin{document}
%
%

\author{Katsunori Kawamura}
\date{xxx,xxx, 2006}



\def\autherp{Katsunori  Kawamura}
\def\emailp{e-mail:kawamura@kurims.kyoto-u.ac.jp.}
\def\addressp{College of Science and Engineering Ritsumeikan University,\\
1-1-1 Noji Higashi, Kusatsu, Shiga 525-8577,Japan
}

%
%
\begin{center}
{\Large \titlep}

\ww
\autherp
\footnote{\emailp}

\noindent
{\it \addressp}
\quad \\

\end{center}

%
%
\begin{abstract}
In our previous articles, we have presented a class of
endomorphisms of the Cuntz algebras which are defined by polynomials of 
canonical generators and their conjugates.
We showed the classification of some case under unitary equivalence 
by help of branching laws of permutative representations. 
In this article, we construct an automaton which is called the Mealy machine
associated with the endomorphism in order to compute its branching law.
We show that the branching law is obtained as outputs
from the machine for the input of information of a given representation.
\end{abstract}

%
%
\sftt{Introduction}
\label{section:first}
In \cite{PE01,PE02}, we introduced a class of endomorphisms of the 
Cuntz algebra $\con$ which are called permutative endomorphisms.
They are given by noncommutative polynomials in canonical generators of $\con$.
Such endomorphisms were motivated by an interest of 
the following endomorphism $\rho_{\nu}$ of $\co{3}$ discovered by Noboru Nakanishi:
%
%
\begin{equation}
\label{eqn:nakeq}
\left\{
\begin{array}{c}
\rho_{\nu}(s_{1})\equiv s_{2}s_{3}s_{1}^{*}+s_{3}s_{1}s_{2}^{*}+s_{1}s_{2}s_{3}^{*},\\
\\
\rho_{\nu}(s_{2})\equiv s_{3}s_{2}s_{1}^{*}+s_{1}s_{3}s_{2}^{*}+s_{2}s_{1}s_{3}^{*},\\
\\
\rho_{\nu}(s_{3})\equiv s_{1}s_{1}s_{1}^{*}+s_{2}s_{2}s_{2}^{*}+s_{3}s_{3}s_{3}^{*}
\end{array}
\right.
\end{equation}
where $s_{1},s_{2},s_{3}$ are canonical generators of $\co{3}$.
Because $\rho_{\nu}(s_{1}),\rho_{\nu}(s_{2}),\rho_{\nu}(s_{3})$
satisfy the relation of canonical generators of $\co{3}$,
we can verify that $\rho_{\nu}$ is an endomorphism of $\co{3}$.
$\rho_{\nu}$ is very concrete but its property is not so clear.
In Theorem 1.2 of \cite{PE01}, 
we proved that $\rho_{\nu}$ is irreducible but not an automorphism
by using branching laws of $\rho_{\nu}$ with respect to permutative representations.
Especially, $\rho_{\nu}$ is not unitarily equivalent
to the canonical endomorphism of $\co{3}$.

In general, representations of C$^{*}$-algebras do not have unique 
decomposition (up to unitary equivalence) into sums or integrals of irreducibles. 
However, the permutative representations of $\con$ do \cite{BJ,DaPi2,DaPi3}.
Because a representation arising from
the right transformation of a permutative representation by 
a permutative endomorphism is also a permutative representation,
their branching laws make sense.
By such branching laws, permutative endomorphisms are characterized 
and classified effectively.

%
%
\begin{defi}
\label{defi:first}
Let $s_{1},\ldots,s_{N}$ be canonical generators of $\con$ and
$({\cal H},\pi)$ be a representation of $\con$.
\begin{enumerate}
\item
$({\cal H},\pi)$ is a {\it permutative representation} of $\con$
if there is a complete orthonormal basis $\{e_{n}\}_{n\in\Lambda}$
of ${\cal H}$ and a family $f=\{f_{i}\}_{i=1}^{N}$
of maps on $\Lambda$ such that $\pi(s_{i})e_{n}=e_{f_{i}(n)}$
for each $n\in\Lambda$ and $i\edot$.
\item
For $J=(j_{i})_{i=1}^{k}\in\nset{k}$,
$({\cal H},\pi)$ is $P(J)$ if there is a unit cyclic vector $\Omega\in {\cal H}$
such that $\pi(s_{J})\Omega=\Omega$ and 
$\{\pi(s_{j_{i}}\cdots s_{j_{k}})\Omega\}_{i=1}^{k}$ is 
an orthonormal family in ${\cal H}$ where $s_{J}\equiv s_{j_{1}}\cdots s_{j_{k}}$.
\item
$({\cal H},\pi)$ is a {\it cycle} if
there is $J\in\nset{k}$ such that $({\cal H},\pi)$ is $P(J)$.
\end{enumerate}
\end{defi}

\noindent
For any $J\in \nset{k}$, $P(J)$ exists uniquely up to unitary equivalence.
In Theorem 1.3 of \cite{PE02},  we showed the following:
%
%
\begin{Thm}
\label{Thm:mainzero}
Let ${\goth S}_{N,l}$ be the set of all permutations on the set $\nset{l}$.
For $\sigma\in {\goth S}_{N,l}$,
let $\psi_{\sigma}$ be the endomorphism of $\con$ defined by
%
%
\begin{equation}
\label{eqn:end}
\psi_{\sigma}(s_{i})\equiv u_{\sigma}s_{i}\quad(i\edot)
\end{equation}
where $u_{\sigma}\equiv \sum_{J\in \nset{l}}s_{\sigma(J)}(s_{J})^{*}$.
If a representation $({\cal H},\pi)$ of $\con$ 
is $P(J)$ for $J\in\nset{k}$ and $\sigma\in {\goth S}_{N,l}$,
then there are $J_{1},\ldots,J_{M}\in\bigcup_{m\geq 1}\nset{m}$ and
subrepresentations $\pi_{1},\ldots,\pi_{M}$ of $\pi\circ \psi_{\sigma}$ such that
%
%
\begin{equation}
\label{eqn:branching}
\pi\circ \psi_{\sigma}=\pi_{1}\oplus\cdots\oplus\pi_{M},
\end{equation}
$\pi_{i}$ is $P(J_{i})$ and $J_{i}\in \coprod_{n=1}^{N^{l-1}}\nset{nk}$ 
for $i=1,\ldots,M$.
Further $1\leq M\leq N^{l-1}$.
\end{Thm}

\noindent
$\psi_{\sigma}$ in (\ref{eqn:end}) is called the {\it permutative endomorphism} 
of $\con$ by $\sigma$.
The canonical endomorphism of $\con$ and $\rho_{\nu}$ in (\ref{eqn:nakeq})
are permutative endomorphisms.

By the uniqueness of decomposition of permutative representation,
the rhs in (\ref{eqn:branching}) is unique up to unitary equivalence.
When $({\cal H},\pi)$ is $P(J)$ and $\rho\in {\rm End}\con$,
we denote $({\cal H},\pi\circ \rho)$ by $P(J)\circ \rho$ simply.
Then (\ref{eqn:branching}) can be rewritten as follows:
%
%
\begin{equation}
\label{eqn:simpdeco}
P(J)\circ \psi_{\sigma}=P(J_{1})\oplus\cdots\oplus P(J_{M}).
\end{equation}
We call (\ref{eqn:simpdeco}) by the {\it branching law} of $\psi_{\sigma}$ 
with respect to $P(J)$.
The  branching law of $\psi_{\sigma}$ is unique up to
unitary equivalence of $\psi_{\sigma}$.
Concrete such branching laws are already given in \cite{PE01,PE02} 
by direct computation.
These branching laws are interesting subjects themselves
and they are useful to classify endomorphisms effectively.
On the other hand, an automaton is a typical object to consider algorithm of
computation in the computer science \cite{Eilenberg,Ginz,hopul,Mealy}.
An automaton is a machine which changes the internal state by an input.
A Mealy machine is a kind of automaton with output.

In this article, we show a better algorithm to compute branching law,
that is, an algorithm to seek $J_{1},\ldots,J_{M}$ from a given $J$ in
(\ref{eqn:simpdeco}) by reducing the problem to a semi-Mealy machine 
$\mms$ as an input ($=J$) and outputs ($=J_{1},\ldots,J_{M}$):

\noindent
\def\chart{
\textsf{{\small
\put(100,100){\framebox(800,300)[c]{$\textsf{M}_{\sigma}$}}
\put(-700,250){\vector(1,0){500}}
\put(-1200,200){{\Large $a_{J}$}}
\put(-1450,400){Input word}
\put(1200,250){\vector(1,0){500}}
\put(2000,200){{\Large $b_{J_{1}},\ldots,b_{J_{M}}$}}
\put(1950,400){Output words}
\put(-90,500){Semi-Mealy machine}
}
}
}
%
%
\setlength{\unitlength}{.022566mm}
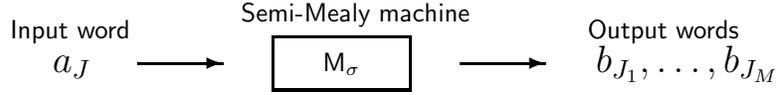
\begin{picture}(1001,650)(0,50)
\thicklines
\put(2000,0){\chart}
\end{picture}

If $J=J_{0}^{r}$, that is, $J$ is a sequence of $r$-times repetition 
of a sequence $J_{0}\in\nset{k^{'}}$ and $r\geq 2$, then there 
are $z_{1},\ldots,z_{r}\in U(1)$ such that
$P(J)=\bigoplus_{j=1}^{r}P(J_{0})\circ \gamma_{z_{j}}$
where $\gamma$ is the gauge action on $\con$ by Theorem 2.4 (iv) in \cite{PE02}.
Because $\gamma_{z}\circ \psi_{\sigma}=\psi_{\sigma}\circ \gamma_{z}$ for each $z$,
the branching law of $P(J)\circ \psi_{\sigma}$
is reduced to that of $P(J_{0})\circ \psi_{\sigma}$.
Therefore it is sufficient to show the case that $J$ is {\it nonperiodic},
that is, $J$ is impossible to be written as $J_{0}^{r}$ for $r\geq 2$.
Hence we assume that $J$ is nonperiodic.

For $\sigma\in {\goth S}_{N,l}$ with $l\geq 2$ and $J\in\nset{l}$,
we define $\sigma_{1}(J),\ldots,\sigma_{l}(J)\in\nset{}$ 
by $\sigma(J)=(\sigma_{1}(J),\ldots,\sigma_{l}(J))$
and let $\sigma_{n,m}(J)\equiv (\sigma_{n}(J),\ldots,
\sigma_{m}(J))$ for $1\leq n<m\leq l$.
Define $\nset{0}\equiv \{0\}$ for convenience.
%
%
\begin{defi}
\label{defi:mealy}
For $\sigma\in {\goth S}_{N,l}$,
a data $\mms\equiv(Q,\Sigma,\Delta,\delta,\lambda)$ 
is called the {\it semi-Mealy machine by $\sigma$} if 
$Q,\Sigma,\Delta$ are finite sets,
\[Q\equiv\{q_{K}:K\in \nset{l-1}\},\quad
\Sigma\equiv\{a_{j}\}_{j=1}^{N},\quad
\Delta\equiv \{b_{j}\}_{j=1}^{N}\]
and two maps $\delta:Q\times\Sigma^{*}\to Q$,
$\lambda:Q\times \Sigma^{*}\to \Delta^{*}$ are defined by
\[
\delta(q_{K},a_{i})\equiv
\left\{
\begin{array}{ll}
q_{0} \quad&(l=1),\\
&\\
q_{(\sigma^{-1})_{2,l}(K,i)}\quad&(l\geq 2),\\
\end{array}
\right.\quad \!\!
\lambda(q_{K},a_{i})\equiv
\left\{
\begin{array}{ll}
b_{\sigma^{-1}(i)} \quad&(l=1),\\
&\\
b_{(\sigma^{-1})_{1}(K,i)}\quad&(l\geq 2)\\
\end{array}
\right.
\]
for $i\edot$ and $K\in\nset{l-1}$ where $\Sigma^{*}$ and
$\Delta^{*}$ are free semigroups generated by $\Sigma$ and $\Delta$, respectively.
\end{defi}

\noindent
We posteriori define $\delta(q,wa)\equiv \delta(\delta(q,w),a)$
and $\lambda(q,wa)\equiv \lambda(q,w)\lambda(\delta(q,w),a)$
for $q\in Q$, $w\in \Sigma^{*}$ and $a\in \Sigma$.
For a given $J=(j_{i})_{i=1}^{k}\in\nset{k}$,
define $Q_{J}\equiv \{q\in Q:\mbox{there exists } n\in{\bf N}\,s.t.\, 
\delta(q,(a_{J})^{n})=q\}$ where 
$a_{J}\equiv a_{j_{1}}\cdots a_{j_{k}}\in \Sigma^{*}$
and define an equivalence relation $\sim$ in $Q_{J}$ by $q\sim q^{'}$
if there is $n\in {\bf N}$ such that $\delta(q,(a_{J})^{n})=q^{'}$.
Define $[q]\equiv\{q^{'}\in Q_{J}:q\sim q^{'}\}$.
Then $[q]$ is a cyclic component of $Q_{J}$ 
with respect to the iteration of the right action of $a_{J}$ by $\delta$.
There are $p_{1},\ldots,p_{M}\in Q_{J}$ such that 
the set $Q_{J}$ of periodic points is decomposed into orbits as follows:
%
%
\begin{equation}
\label{eqn:qdeco}
Q_{J}=[p_{1}]\sqcup\cdots \sqcup [p_{M}].
\end{equation}
Under these preparations, the main theorem is given as follows:
%
%
\begin{Thm}
\label{Thm:mainone}
If $J$ is nonperiodic,
then $J_{1},\ldots,J_{M}$ in (\ref{eqn:simpdeco}) are obtained by 
\[b_{J_{i}}=\lambda(p_{i},(a_{J})^{r_{i}}) \quad (i=1,\ldots,M)\]
where $p_{1},\ldots,p_{M}\in Q_{J}$ are taken as (\ref{eqn:qdeco}) 
and $r_{i}\equiv \#[p_{i}]$ for $i=1,\ldots,M$.
\end{Thm}

\noindent
In Theorem \ref{Thm:mainone}, if $p_{1}^{'},\ldots,p_{M}^{'}$
satisfy (\ref{eqn:simpdeco}) and $[p_{i}^{'}]=[p_{i}]$ for each $i$,
then the associated $J_{1}^{'},\ldots,J_{M}^{'}$
satisfy that $P(J_{i}^{'})=P(J_{i})$ for each $i$.
We show a more practical algorithm to compute branching laws 
by using the Mealy diagram as follows:

The {\it transition diagram (Mealy diagram)} 
${\cal D}(\mathsf{M})$ of a semi-Mealy machine $\mathsf{M}
=(Q,\Sigma,\Delta,\delta,\lambda)$
is a directed graph with labeled edges, which has a set $Q$ of vertices and a set 
$E\equiv \{(q,\delta(q,a),a)\in Q\times Q\times \Sigma:q\in Q,\, a\in \Sigma\}$
of directed edges with labels.
The meaning of $(q,\delta(q,a),a)$ is an edge from $q$ to $\delta(q,a)$
with a label ``$a/\lambda(q,a)$" for $a\in\Sigma$:

\noindent
%
%
\def\jei{\put(0,0){\oval(600,300)}\put(-70,-40){$q$}}
\def\jeitwo{\put(0,0){\oval(600,300)}\put(-70,-50){$p$}}
\def\figures{
\put(0,0){\jei}
\put(1500,0){\jeitwo}
\put(300,0){\vector(1,0){900}}
\put(600,30){$^{a/b}$}
}
%
%
\setlength{\unitlength}{.022566mm}
\begin{picture}(1001,450)(99,-250)
\thicklines
\put(300,-50){$\delta(q,a)=p$, \quad $\lambda(q,a)=b$ \quad
$\Longleftrightarrow$}
\put(3200,0){\figures}
\end{picture}

For $\rho_{\nu}$ in (\ref{eqn:nakeq}), 
we compute branching laws by the semi-Mealy machine.
Define $\sigma_{0}\in{\goth S}_{3,2}$ by
${\small \begin{array}{c|ccccccccc}
J&11&12&13&21&22&23&31&32&33\\ \hline
\sigma_{0}(J)&23&31&12&32&13&21&11&22&33\\
\end{array}}$.
Then $\rho_{\nu}= \psi_{\sigma_{0}}$ and
$\mathsf{M}_{\sigma_{0}}=(\{q_{1},q_{2},q_{3}\},
\{a_{1},a_{2},a_{3}\}$, $\{b_{1},b_{2},b_{3}\},
\delta,\lambda)$ is given as follows:

{\small
\[
\begin{array}{c|c|c|c|c|c|c}
p&\delta(p,a_{1})&\delta(p,a_{2})&\delta(p,a_{3})
&\lambda(p,a_{1})&\lambda(p,a_{2})&\lambda(p,a_{3})\\
\hline
q_{1}&q_{1}&q_{3}&q_{2}&b_{3}&b_{1}&b_{2}\\
q_{2}&q_{3}&q_{2}&q_{1}&b_{2}&b_{3}&b_{1}\\
q_{3}&q_{2}&q_{1}&q_{3}&b_{2}&b_{1}&b_{3}\\
\end{array}
\]
}

\noindent
From this, ${\cal D}(\mathsf{M}_{\sigma_{0}})$ is as follows:

\noindent
%
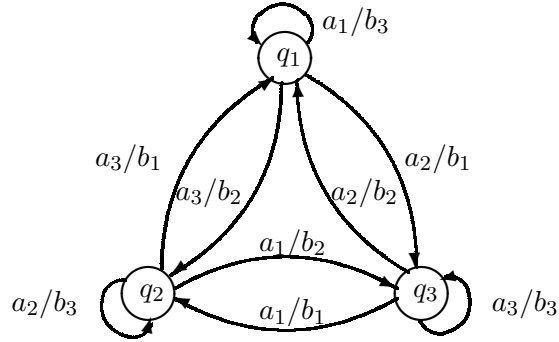
\def\bcircle#1{\circle{300}\put(-50,-20){$q_{#1}$}}
\def\labee#1#2{\put(0,0){$a_{#1}/b_{#2}$}}
\def\numbers{
\put(-800,0){\bcircle{2}}
\put(0,1380){\bcircle{1}}
\put(800,0){\bcircle{3}}
}
\def\edgethree{
\thicklines
\qbezier(-650,-30)(0,-400)(650,-30)
\qbezier(-650,30)(0,400)(650,30)
\put(-650,-20){\vector(-2,1){0}}
\put(650,30){\vector(2,-1){0}}
}
\def\edges{
\thicklines
\put(0,0){\edgethree}
\put(230,700){\rotatebox{120}{\edgethree}}
\put(-560,690){\rotatebox{-120}{\edgethree}}
}
\def\semi{
\qbezier(150,0)(0,150)(-150,0)
\qbezier(150,0)(200,-75)(150,-150)
\qbezier(-150,0)(-200,-75)(-150,-150)
\put(-110,-170){\vector(1,-1){0}}}
\def\labels{
\put(220,1550){\labee{1}{3}}
\put(-1600,-100){\labee{2}{3}}
\put(1210,-100){\labee{3}{3}}
%
\put(-1110,750){\labee{3}{1}}
\put(700,750){\labee{2}{1}}
\put(-150,-150){\labee{1}{1}}
\put(-150,250){\labee{1}{2}}
\put(-640,550){\labee{3}{2}}
\put(270,550){\labee{2}{2}}
}
%
\setlength{\unitlength}{.022713mm}
\begin{picture}(3011,2120)(-2299,-360)
\thicklines
\put(0,0){\labels}
\put(0,0){\numbers}
\put(0,0){\edges}
%
\put(-20,1620){\semi}
\put(850,-100){\rotatebox{-120}{\semi}}
\put(-1055,-130){\rotatebox{120}{\semi}}
\end{picture}

\noindent
According to Theorem \ref{Thm:mainone},
we compute branching laws for $\rho_{\nu}$
by ${\cal D}(\mathsf{M}_{\sigma_{0}})$.
When the input word is  $a_{1}$,
$\delta(q_{1},a_{1})=q_{1}$,
$\delta(q_{2},a_{1})=q_{3}$,
$\delta(q_{3},a_{1})=q_{2}$.
Therefore $Q_{1}=[q_{1}]\sqcup[q_{2}]$, $r_{1}=1$, $r_{2}=2$ and
there are two cycles $q_{1}$ and $q_{2}q_{3}$ in $Q$ with respect to $a_{1}$.
From this, we have output words,
$\lambda(q_{1},a_{1})=b_{3}$ and $\lambda(q_{2},(a_{1})^{2})=b_{2}b_{1}$.
Hence $P(1)\circ \rho_{\nu}=P(3)\oplus P(21)=P(3)\oplus P(12)$.
where we use a fact that
$P(j_{p(1)},\ldots,j_{p(k)})=P(j_{1},\ldots,j_{k})$ for each $p\in {\bf Z}_{k}$.
Further the following holds:

{\footnotesize
\[
\begin{array}{c|c|c|c}
\mbox{input}&\mbox{cycles}&\mbox{outputs}
&\mbox{branching law}\\
\hline
a_{1}&q_{1},q_{2}q_{3}
&b_{3},b_{2}b_{1}&P(1)\circ \rho_{\nu}=P(3)\oplus P(12)\\
\hline
a_{1}a_{2}&q_{1}q_{1}q_{3}q_{2}q_{2}q_{3}
&b_{3}b_{1}b_{1}b_{3}b_{2}b_{2}
&P(12)\circ \rho_{\nu}=P(113223)\\
\hline
a_{1}a_{2}a_{3}&q_{1}q_{1}q_{3}q_{3}q_{2}q_{2},q_{2}q_{3}q_{1}
&b_{3}b_{1}b_{3}b_{1}b_{3}b_{1},b_{2}b_{2}b_{2}
&P(123)\circ \rho_{\nu}=P(131313)\oplus P(222)\\
\hline
a_{1}a_{3}a_{2}&q_{1}q_{1}q_{2}q_{2}q_{3}q_{3},q_{3}q_{2}q_{1}
&b_{3}b_{2}b_{3}b_{2}b_{3}b_{2},b_{1}b_{1}b_{1}
&P(132)\circ \rho_{\nu}=P(232323)\oplus P(111)\\
\end{array}
\]
}

In $\S$\ref{section:third}, we rewrite branching laws by
branching function systems and their transformations, and
we review known facts about endomorphisms.
$\S$\ref{section:fourth} is devoted to prove Theorem \ref{Thm:mainone}
by branching function systems. 
In $\S$\ref{section:fifth}, we show examples of Mealy diagram of 
the semi-Mealy machine $\mms$ and branching laws of $\psi_{\sigma}$
for concrete $\sigma\in {\goth S}_{N,l}$.

%
%
\sftt{Branching function systems}
\label{section:third}
In order to compute branching laws of endomorphisms,
we introduce branching function systems and their transformations by permutations.

Let $\nset{*}_{1}\equiv \bigcup_{k\geq 1}\nset{k}$.
For $J\in \nset{*}_{1}$, the {\it length} of $J$ is defined by
$k$ when $J\in \nset{k}$.
For $J_{1}=(j_{1},\ldots,j_{k}),J_{2}=(j_{1}^{'},\ldots,j_{l}^{'})$,
let $J_{1}\cup J_{2}\equiv(j_{1},\ldots,j_{k},j_{1}^{'},\ldots,j_{l}^{'})$.
Especially, we define $(i,J)\equiv (i)\cup J$ for convenience.
For $J$ and $k\geq 2$, $J^{k}=J\cup\cdots\cup J$ ($k$-times).
For $J=(j_{1},\ldots,j_{k})$ and $\tau\in {\bf Z}_{k}$,
define $\tau(J)\equiv (j_{\tau(1)}, \ldots,j_{\tau(k)})$.
For $J_{1},J_{2}\in \nset{*}_{1}$,
$J_{1}\sim J_{2}$ if there are $k\geq 1$ and $\tau\in {\bf Z}_{k}$
such that $J_{1},J_{2}\in\nset{k}$ and $\tau(J_{1})=J_{2}$.
For $J_{1}=(j_{1},\ldots,j_{k}),J_{2}=(j_{1}^{'},\ldots,j_{k}^{'})$, 
$J_{1}\prec J_{2}$ if $\sum_{l=1}^{k}(j_{l}^{'}-j_{l})N^{k-l}\geq 0$.
$J\in\nset{*}_{1}$ is {\it minimal} if $J\prec J^{'}$ for each 
$J^{'}\in\nset{*}_{1}$ such that $J\sim J^{'}$.
Define $[1,\ldots,N]^{*}\equiv \{J\in\nset{*}_{1}:
J\mbox{ is minimal and nonperiodic}\}$.
$[1,\ldots,N]^{*}$ is in one-to-one correspondence
with the set of all equivalence classes of 
nonperiodic elements in $\nset{*}_{1}$
with respect to the equivalence relation $\sim$.

Let $\Lambda$ be an infinite set and $N\geq 2$.
$f=\{f_{i}\}_{i=1}^{N}$ is a {\it branching function system on $\Lambda$}
if $f_{i}$ is an injective transformation
on $\Lambda$ for $i\edot$ such that
a family of their images coincides a partition of $\Lambda$.
Let $\bfsnl$ be the set of all branching function systems on $\Lambda$.
$f=\{f_{i}\}_{i=1}^{N}\in{\rm BFS}_{N}(\Lambda_{1})$ 
and $g=\{g_{i}\}_{i=1}^{N}\in{\rm BFS}_{N}(\Lambda_{2})$ 
are {\it equivalent} if there is a bijection
$\varphi$ from $\Lambda_{1}$ to $\Lambda_{2}$ such that
$\varphi\circ f_{i}\circ \varphi^{-1}=g_{i}$ for $i\edot$.
For $f=\{f_{i}\}_{i=1}^{N}$,
we denote $f_{J}\equiv f_{j_{1}}\circ\cdots\circ f_{j_{k}}$
when $J=(j_{1},\ldots,j_{k})\in \nset{k}$ and define $f_{0}\equiv id$.
For $x,y\in \Lambda$, $x\sim y$ (with respect to $f$)
if there are $J_{1},J_{2}\in\nset{*}$ and $z\in \Lambda$
such that $f_{J_{1}}(z)=x$ and $f_{J_{2}}(z)=y$.
For $x\in \Lambda$, define $A_{f}(x)\equiv \{y\in \Lambda:x\sim y\}$.
$f=\{f_{i}\}_{i=1}^{N}\in\bfsnl$ is {\it cyclic} 
if there is an element $x\in \Lambda$ such that $\Lambda=A_{f}(x)$.
$\{n_{1},\ldots,n_{k}\}\subset \Lambda$ is a {\it cycle} of $f$
if there is $J=(j_{1},\ldots,j_{k})$ such that 
$f_{j_{1}}(n_{1})=n_{k},f_{j_{2}}(n_{2})=n_{1},\ldots,f_{j_{k}}(n_{k})=n_{k-1}$.
$f$ has a {\it cycle} if there is a cycle of $f$ in $\Lambda$.

Let $\Xi$ be a set and $\Lambda_{\omega}$ be an infinite set for $\omega\in\Xi$.
For $f^{[\omega]}=\{f_{i}^{[\omega]}\}_{i=1}^{N}\in{\rm BFS}_{N}(\Lambda_{\omega})$,
$f$ is the {\it direct sum} of $\{f^{[\omega]}\}_{\omega\in\Xi}$
if $f=\{f_{i}\}_{i=1}^{N}\in {\rm BFS}_{N}(\Lambda)$
for a set $\Lambda\equiv \coprod_{\omega\in\Xi}\Lambda_{\omega}$ which
is defined by $f_{i}(n)\equiv f^{[\omega]}_{i}(n)$ when $n\in\Lambda_{\omega}$
for $i\edot$ and $\omega\in\Xi$.
For $f\in\bfsnl$,
$f=\bigoplus_{\omega\in\Xi}f^{[\omega]}$ is a {\it decomposition} of $f$
into a family $\{f^{[\omega]}\}_{\omega\in\Xi}$
if there is a family $\{\Lambda_{\omega}\}_{\omega\in\Xi}$ 
of subsets of $\Lambda$ such that $f$ is the direct sum of 
$\{f^{[\omega]}\}_{\omega\in\Xi}$.
For each $f=\{f_{i}\}_{i=1}^{N}\in\bfsnl$,
there is a decomposition $\Lambda=\coprod_{\omega\in\Xi}\Lambda_{\omega}$
such that $\#\Lambda_{\omega}=\infty$,
$f|_{\Lambda_{\omega}}\equiv \{f_{i}|_{\Lambda_{\omega}}\}_{i=1}^{N}
\in{\rm BFS}_{N}(\Lambda_{\omega})$ and $f|_{\Lambda_{\omega}}$
is cyclic for each $\omega\in\Xi$.
%
%
\begin{defi}
\label{defi:pj}
\begin{enumerate}
\item
For $J\in \nset{k}$, $f\in\bfsnl$ is $P(J)$
if $f$ is cyclic and has a cycle $\{n_{1},\ldots,n_{k}\}$
such that $f_{J}(n_{k})=n_{k}$.
\item
For $f\in\bfsnl$ and $J\in\nset{*}_{1}$,
$g$ is a $P(J)$-component of $f$ 
if $g$ is a direct sum component of $f$ and $g$ is $P(J)$.
\end{enumerate}
\end{defi}

\noindent
For $f\in\bfsnl$ and $\Lambda_{1},\Lambda_{2}\subset \Lambda$,
if $f|_{\Lambda_{i}}$ is $P(J_{i})$ for $i=1,2$,
then either $\Lambda_{1}\cap \Lambda_{2}=\emptyset$ or $\Lambda_{1}=\Lambda_{2}$.

Recall ${\goth S}_{N,l}$ in Theorem \ref{Thm:mainzero}.
For $\sigma\in {\goth S}_{N,l}$ and $f=\{f_{i}\}_{i=1}^{N}\in\bfsnl$,
define $f^{(\sigma)}=\{f^{(\sigma)}_{i}\}_{i=1}^{N}\in\bfsnl$ by
%
%
\begin{equation}
\label{eqn:sigmaf}
f^{(\sigma)}_{i}\equiv f_{\sigma(i)}\quad(l=1),\quad
f^{(\sigma)}_{i}(f_{J}(n))\equiv f_{\sigma(i,J)}(n)\quad(l\geq 2)
\end{equation}
for $n\in\Lambda$, $i\edot$ and $J\in \nset{l-1}$.
If $\sigma\in {\goth S}_{N}={\goth S}_{N,1}$ and $f\in\bfsnl$ is $P(J)$,
then $f^{(\sigma)}$ is $P(J_{\sigma^{-1}})$ where 
$J_{\sigma^{-1}}\equiv\left(\sigma^{-1}(j_{1}),\ldots,\sigma^{-1}(j_{k})\right)$
for $J=(j_{1},\ldots,j_{k})$.
For any $J\in\nset{*}_{1}$, 
there is $f\in\bfsnl$ for some set $\Lambda$ such that $f$ is $P(J)$.
In this case, for $\sigma\in {\goth S}_{N,l}$,
there is $1\leq M\leq N^{l-1}$ such that
$f^{(\sigma)}$ is decomposed into a direct sum of $M$ cycles
by Lemma 2.2 in \cite{PE02}.
Furthermore, the length of each cycle is a multiple of that of $J$.

For $N\geq 2$, let $\con$ be the {\it Cuntz algebra} \cite{C}, that is,
the C$^{*}$-algebra which is universally generated by
$s_{1},\ldots,s_{N}$ satisfying $s^{*}_{i} s_j=\delta_{ij}I$ for 
$i,j\edot$ and $s_1 s^{*}_1+\cdots+s_N s^{*}_N=I$.
In this article, any representation and endomorphism are assumed unital and $*$-preserving.

$(l_{2}(\Lambda),\pi_{f})$ is the {\it permutative representation
of $\con$ by $f=\{f_{i}\}_{i=1}^{N}\in\bfsnl$}
if $\pi_{f}(s_{i})e_{n}\equiv e_{f_{i}(n)}$
for $n\in\Lambda$ and $i\edot$.
For $J\in \nset{*}_{1}$, $P(J)$ in Definition \ref{defi:first}
is irreducible if and only if $J$ is nonperiodic.
For $J_{1},J_{2}\in \nset{*}_{1}$, $P(J_{1})\sim P(J_{2})$
if and only if $J_{1}\sim J_{2}$ where $P(J_{1})\sim P(J_{2})$ means
the unitary equivalence of two representations which satisfy
the condition $P(J_{1})$ and $P(J_{2})$, respectively.
$[1,\ldots,N]^{*}$ is in one-to-one correspondence
with the set of equivalence classes of irreducible 
permutative representations of $\con$ with a cycle.
If $f\in\bfsnl$ and $g\in{\rm BFS}_{N}(\Lambda^{'})$ satisfy $f\sim g$,
then $(l_{2}(\Lambda),\pi_{f})\sim(l_{2}(\Lambda^{'}),\pi_{g})$.
If $f$ is cyclic, then $(l_{2}(\Lambda),\pi_{f})$ is cyclic.
If $f$ is $P(J)$, then $(l_{2}(\Lambda),\pi_{f})$ is $P(J)$.
If $\Lambda=\Lambda_{1}\sqcup \Lambda_{2}$ and 
$f^{(i)}\equiv f|_{\Lambda_{i}}\in{\rm BFS}_{N}(\Lambda_{i})$ for $i=1,2$,
then $(l_{2}(\Lambda),\pi_{f})\sim
(l_{2}(\Lambda_{1}),\pi_{f^{(1)}})\oplus(l_{2}(\Lambda_{2}),\pi_{f^{(2)}})$.

Let $\enda$ be the set of all unital $*$-endomorphisms of 
a unital $*$-algebra ${\cal A}$. 
For $\rho\in\enda$, $\rho$ is {\it proper} if $\rho({\cal A})\ne {\cal A}$.
$\rho$ is {\it irreducible} if $\rho({\cal A})^{'}\cap  {\cal A}={\bf C}I$
where $\rho({\cal A})^{'}\cap  {\cal A}
\equiv \{x\in {\cal A}:\mbox{for all } a\in{\cal A},\,\rho(a)x=x\rho(a)\}$.
$\rho$ and $\rho^{'}$ are {\it equivalent}
if there is a unitary $u\in{\cal A}$
such that $\rho^{'}={\rm Ad}u\circ \rho$.
In this case, we denote $\rho\sim\rho^{'}$.
Let ${\rm Rep}{\cal A}$ ({\it resp.} ${\rm IrrRep}{\cal A}$) be the set of all
unital ({\it resp.} irreducible) $*$-representations of ${\cal A}$.
We simply denote $\pi$ for $({\cal H},\pi) \in{\rm Rep}{\cal A}$.
If $\rho,\rho^{'}\in {\rm End}{\cal A}$ and $\pi,\pi^{'}\in {\rm Rep}{\cal A}$
satisfy $\rho\sim\rho^{'}$ and $\pi\sim \pi^{'}$,
then $\pi\circ \rho\sim  \pi^{'}\circ \rho^{'}$.
Assume that ${\cal A}$ is simple.
If there is $\pi\in {\rm IrrRep}{\cal A}$ such that
$\pi\circ \rho\in {\rm IrrRep}{\cal A}$, then $\rho$ is irreducible.
If there is $\pi\in {\rm Rep}{\cal A}$ such that
$\pi\circ\rho\not\sim \pi\circ\rho^{'}$, then $\rho\not\sim\rho^{'}$.
If there is $\pi\in {\rm IrrRep}{\cal A}$
such that $\pi\circ\rho\not\in {\rm IrrRep}{\cal A}$, then $\rho$ is proper.

For $\psi_{\sigma}$ in (\ref{eqn:end}), define
%
%
\begin{equation}
\label{eqn:edef}
E_{N,l}\equiv\{\psi_{\sigma}\in\endcon:\sigma\in {\goth S}_{N,l}\}\quad(l\geq 1).
\end{equation}
If $\sigma\in {\goth S}_{N}$,
then $\psi_{\sigma}$ is an automorphism of $\con$ which satisfies
$\psi_{\sigma}(s_{i})=s_{\sigma(i)}$ for $i\edot$.
Especially, if $\sigma=id$, then $\psi_{id}=id$.
If $\sigma\in {\goth S}_{N,2}$ is defined by
$\sigma(i,j)\equiv (j,i)$ for $i,j\edot$,
then $\psi_{\sigma}$ is just the canonical endomorphism of $\con$.
For $\sigma\in {\goth S}_{N,l}$ and $f\in\bfsnl$,
$\pi_{f}\circ \psi_{\sigma}=\pi_{f^{(\sigma)}}$
where $f^{(\sigma)}$ is in (\ref{eqn:sigmaf}).
If $\rho$ is a permutative endomorphism and 
$({\cal H},\pi)$ is a permutative representation of $\con$,
then $\pi\circ\rho$ is also a permutative representation.

A representation $({\cal H},\pi)$ of $\con$ has a {\it $P(J)$-component}
if $({\cal H},\pi)$ has a subrepresentation
$({\cal H}_{0},\pi|_{{\cal H}_{0}})$ which is $P(J)$.
A component of a representation $P(J)\circ \rho$
of $\con$ means a subrepresentation of $({\cal H},\pi)$
which is equivalent to $P(J^{'})$ for some $J^{'}$.

For comparison of the method to find $(J_{i})_{i=1}^{M}$ in (\ref{eqn:simpdeco})
for a given $J$,
we show the usual method to determine $(J_{i})_{i=1}^{M}$ as follows:
(a) Prepare a representation $({\cal H},\pi)$ which is $P(J)$.
We often take ${\cal H}=\ltn$ and $\pi=\pi_{f}$ 
for suitable branching function system $f$ on ${\bf N}$.
(b) Compute $\pi(\psi_{\sigma}(s_{i}))e_{n}$ for each $n\in {\bf N}$ and $i\edot$.
By the proof of Lemma 2.2 in \cite{PE02}, we see that it is sufficient to
check for $1\leq n\leq N^{l-1}k$ when $|J|=k$.
(c) Find all cycles in ${\cal H}$ by using results in (b).
In this way, the direct computation of branching law is too much of a bother
because of a great number of calculated amount when $N,k,l$ are large.

%
%
\sftt{Proof of Theorem \ref{Thm:mainone}}
\label{section:fourth}
In this section, we assume that $\sigma\in{\goth S}_{N,l}$,
$l\geq 2$, $J=(j_{i})_{i=1}^{k}\in\nset{k}$ and $J$ is nonperiodic.
For $r\geq 2$, extend $J=(j_{i})_{i=1}^{k}$ as $(j_{n})_{n=1}^{r\cdot k}$ by
$j_{k(c-1)+i}\equiv j_{i}$ for each $c=1,\ldots,r$ and $i=1,\ldots,k$
for convenience.
%
%
\begin{lem}
\label{lem:direct}
Let $f\in \bfsnl$ be $P(J)$, $f^{(\sigma)}$ be in (\ref{eqn:sigmaf})
and let $M_{\sigma}=(Q,\Sigma,\Delta,\delta,\lambda)$
be in Definition \ref{defi:mealy}.
For $p\in Q_{J}$, define $r_{J}(p)\in{\bf N}$ by $r_{J}(p)\equiv \#[p]$.
\begin{enumerate}
\item
For $p\in Q_{J}$ and $\alpha\equiv r_{J}(p)\cdot k$,
define $p_{1},\ldots,p_{\alpha}\in Q$
and $T=(t_{i})_{i=1}^{\alpha}\in\nset{\alpha}$ by
$p_{1}\equiv p$, 
$b_{t_{1}}=\lambda(p_{\alpha},a_{j_{\alpha}})$ and 
\[ p_{i}\equiv \delta(p_{i-1},a_{j_{i-1}}),\quad
b_{t_{i}}=\lambda(p_{i-1},a_{j_{i-1}})\quad(i=2,\ldots,\alpha),\]
then 
there is $\Lambda(p)\subset \Lambda$ such that $f^{(\sigma)}|_{\Lambda(p)}$ 
is $P(T)$.
\item
In (i), define $T^{'}\in\nset{\alpha}$ by 
$b_{T^{'}}=\lambda(p,a_{J}^{r_{J}(p)})$.
Then $f^{(\sigma)}|_{\Lambda(p)}$ is $P(T^{'})$.
\item
If there is $\Lambda_{0}$ such that $f^{(\sigma)}|_{\Lambda_{0}}$ is $P(T)$
for $T=(t_{i})_{i=1}^{\alpha}\in\nset{\alpha}$, 
then there is $p\in Q_{J}$ such that $\Lambda_{0}$ is equal to $\Lambda(p)$ in (i).
\item
In (i), $p\sim p^{'}$ if and only if $\Lambda(p)=\Lambda(p^{'})$.
\item
Choose $p_{1},\ldots,p_{M}$ as (\ref{eqn:qdeco}).
Then the decomposition $f^{(\sigma)}=f^{[1]}\oplus\cdots\oplus f^{[M]}$ 
holds as a branching function system 
where $f^{[i]}\equiv f^{(\sigma)}|_{\Lambda(p_{i})}$ for each $i$.
\end{enumerate}
\end{lem}
%
%
\pr
Let $n_{0}\in\Lambda$ such that $f_{J}(n_{0})=n_{0}$.
Because $J$ is nonperiodic, such $n_{0}$ is unique in $\Lambda$.

\noindent
(i)
Let $r\equiv r_{J}(p)$.
There is a sequence $(I_{1},\ldots,I_{\alpha})$ in $\nset{l-1}$ such that
$p_{i}=q_{I_{i}}$ for each $i$. 
By definition of $\delta$ and $\lambda$ and assumption,
%
%
\begin{equation}
\label{eqn:its}
\sigma(t_{1},I_{1})=(I_{\alpha},j_{\alpha}),\sigma(t_{2},I_{2})=(I_{1},j_{1}), 
\ldots,\sigma(t_{\alpha},I_{\alpha})=(I_{\alpha-1},j_{\alpha-1}).
\end{equation}
Define $m(p)\equiv f_{\sigma(t_{1},I)}(n_{0})\in\Lambda$.
Then $m(p)=f_{I_{\alpha}}(f_{j_{\alpha}}(n_{0}))$.
By this and definition of $f^{(\sigma)}$,
we can verify that $f^{(\sigma)}_{T}(m(p))=m(p)$.
Define 
\[m_{\alpha}\equiv m(p),\quad m_{\alpha-1}\equiv f^{(\sigma)}_{t_{\alpha}}(m(p)),
\ldots,m_{1}\equiv f^{(\sigma)}_{(t_{1},\ldots,t_{\alpha})}(m(p))\]
and $\Lambda(p)\equiv \{f^{(\sigma)}_{K}(m(p)):K\in\nset{*}_{1}\}$.
It is sufficient to show that $m_{i}\ne m_{j}$ when $i\ne j$.
By definition,
\[m_{i}=f^{(\sigma)}_{t_{i+1}}(m_{i+1})
=f_{(I_{i},j_{i})}(f_{(j_{i+1},\ldots,j_{\alpha})}(n_{0}))
\quad(i=1,\ldots,\alpha-1)\quad m_{\alpha}=f^{(\sigma)}_{t_{1}}(m_{1}).\]
Assume that $m_{i}=m_{i^{'}}$ and $c\equiv i^{'}-i\geq 0$.
This implies that $m_{\tau(i)}=m_{\tau(i^{'})}$ for each $\tau\in{\bf Z}_{\alpha}$.
From this,
$(I_{\tau(i)},j_{\tau(i)})=(I_{\tau(i^{'})},j_{\tau(i^{'})})$ and
$f^{(\sigma)}_{(t_{i+1},\ldots,t_{\alpha})}(m(p))
=f^{(\sigma)}_{(t_{i^{'}+1},\ldots,t_{\alpha})}(m(p))$.
This implies that $f_{(I_{\alpha},j_{\alpha})}(n_{0})
=f_{(I_{c},j_{c})}(f_{(j_{c+1},\ldots,j_{\alpha})}(n_{0}))$.
Therefore $n_{0}=f_{(j_{c+1},\ldots,j_{\alpha})}(n_{0})$.
By the uniqueness of the cycle in $\Lambda$ with respect to $f$,
$c=k(d-1)$ for $1\leq d\leq r$.
Hence $I_{\tau(i)}=I_{\tau(i+k(d-1))}$ for each $\tau$.
Therefore
$p_{\tau(i)}=q_{I_{\tau(i)}}=q_{I_{\tau(i+k(d-1))}}=p_{\tau(i+k(d-1))}$ 
for each $\tau$.
By the choice of $r$, $d=1$ and $i=i^{'}$. Hence the statement holds.

\noindent
(ii)
We see that $t_{1}^{'}=t_{\alpha},
t_{2}^{'}=t_{1},\ldots,t_{\alpha}^{'}=t_{\alpha-1}$.
Hence $P(T)\sim P(T^{'})$ by definition.

\noindent
(iii)
Fix $\tau\in {\bf Z}_{\alpha}$.
Define $T^{'}=(t_{i}^{'})^{\alpha}_{i=1}\in\nset{\alpha}$ by 
%
%
\begin{equation}
\label{eqn:newt}
t_{i}^{'}\equiv t_{\tau^{-1}(i)}\quad(i=1,\ldots,\alpha).
\end{equation}
Then $f^{(\sigma)}|_{\Lambda_{0}}$ is also $P(T^{'})$
and there is $m_{0}\in\Lambda_{0}$ such that
$f_{T^{'}}^{(\sigma)}(m_{0})=m_{0}$.
Define $m_{\alpha}\equiv m_{0}$ and
$m_{i}\equiv f^{(\sigma)}{(t_{i+1},\ldots,t_{\alpha})}(m_{0})$ for 
$i=1,\ldots,\alpha-1$.
Then $m_{i}\ne m_{i^{'}}$ when $i\ne i^{'}$.
By definition of $f$,
there are $n^{'}\in\Lambda$, $I_{0}\in\nset{l-1}$ and $u_{0}\in\nset{}$
such that $m_{\alpha}=f_{(I_{0},u_{0})}(n^{'})$.
Define a sequence $(I_{i}^{'})_{i=1}^{\alpha}$ in $\nset{l-1}$ and 
$U=(u_{i})_{i=1}^{\alpha}\in\nset{\alpha}$ by
\[I_{\alpha}^{'}\equiv I_{0},\quad u_{\alpha}\equiv u_{0},\quad
(I_{i}^{'},u_{i})\equiv \sigma(t_{i+1}^{'},I_{i+1}^{'})\quad
(i=\alpha-1,\alpha-2,\ldots,1).\]
By assumption, we see that $f_{(I_{\alpha}^{'},u_{\alpha})}(n^{'})
=f_{\sigma(t_{1}^{'},I_{1}^{'})}(f_{U}(n^{'}))$.
By definition of $f$,
$(I_{\alpha}^{'},u_{\alpha})=\sigma(t_{1}^{'},I_{1}^{'})$ and $n^{'}=f_{U}(n^{'})$.
By the uniqueness of cycle in $\Lambda$ with respect to $f$, $U\sim J^{r}$.
Hence there is $\tau^{'}\in{\bf Z}_{\alpha}$ such that
$j_{i}=u_{\tau^{'}(i)}$ for $i=1,\ldots,\alpha$.
Here choose $\tau$ in (\ref{eqn:newt}) by $\tau\equiv \tau^{'}$ and
define $I_{i}\equiv I_{\tau(i)}^{'}$ for each $i$.
Then (\ref{eqn:its}) holds.
From this, we can verify that $p\equiv q_{I_{1}}$ belongs to $Q_{J}$.
Define $m(p)\equiv f_{(I_{\alpha},j_{\alpha})}(n_{0})$ as (i).
Then $n_{0}= f_{(j_{1},\ldots,j_{\tau^{-1}(\alpha)})}(n^{'})$ and
$m_{\alpha}=f^{(\sigma)}_{(t_{\tau^{-1}(1)},\ldots,t_{\alpha})}(m(p))$.
Therefore $m_{\alpha}\in\Lambda(p)$.
Since $m_{\alpha}\in \Lambda_{0}\cap\Lambda(p)$, $\Lambda_{0}=\Lambda(p)$.

\noindent
(iv)
If $p\sim p^{'}$,
then there is $c$ such that $p^{'}=p_{kc+1}$ in (i) and we can verify that
$m(p^{'})=f^{(\sigma)}_{(t_{1+kc},\ldots,t_{\alpha})}(m(p))\in\Lambda(p)$.
Since $m(p^{'})\in \Lambda(p^{'})\cap\Lambda(p)$, $\Lambda(p^{'})=\Lambda(p)$.
 
Assume that $\Lambda(p)=\Lambda(p^{'})$.
Let $m(p),m(p^{'})\in \Lambda$ be in the proof of (i).
Then there are $T,T^{'}\in\nset{*}_{1}$ such that
$f^{(\sigma)}_{T}(m(p))=m(p)$ and $f^{(\sigma)}_{T^{'}}(m(p^{'}))=m(p^{'})$.
Then $f^{(\sigma)}|_{\Lambda(p)}$ is $P(T)$ and
$f^{(\sigma)}|_{\Lambda(p^{'})}$ is $P(T^{'})$.
Since $f^{(\sigma)}|_{\Lambda(p)}=f^{(\sigma)}|_{\Lambda(p^{'})}$, $T^{'}\sim T$.
Assume that $T=(t_{i})_{i=1}^{\alpha}$ and $T^{'}=(t_{i}^{'})_{i=1}^{\alpha}$.
Let $\{m_{i}\}_{i=1}^{\alpha}$ be the cycle 
in $\Lambda(p)$ of $f^{(\sigma)}$ in (i).
By the uniqueness of the cycle in $\Lambda(p)$ with respect to $f^{(\sigma)}$,
$\{m_{i}\}_{i=1}^{\alpha}$ is also the cycle in $\Lambda(p^{'})$ of $f^{(\sigma)}$.
By the proof of (i), $m(p^{'})\in\{m_{i}\}_{i=1}^{\alpha}$.
Hence there is $\tau\in {\bf Z}_{\alpha}$ such that $m(p^{'})=m_{\tau(\alpha)}$.
From this, $t_{i}^{'}=t_{\tau(i)}$ for $i=1,\ldots,\alpha$.
Because $T\sim T^{'}$, $r_{J}(p^{'})=r_{J}(p)$.
Let $r\equiv r_{J}(p)$. 
Assume that $p=q_{I_{1}}$ and $p^{'}=q_{I_{1}^{'}}$.
By definition of $m(p)$ and $m(p^{'})$ and their relation,
we see that $I_{1}^{'}=I_{\tau(1)}$.
Therefore $p^{'}=q_{I_{\tau(1)}}$.
By choice of $p$ and $p^{'}$,
$\delta(p,a_{J}^{r})=p$ and $\delta(p^{'},a_{J}^{r})=p^{'}$.
Because $J$ is nonperiodic, $\tau(i)=i+kc$ for a certain $c$ modulo $\alpha$.
Therefore $p^{'}=q_{I_{\tau(1)}}=q_{I_{1+kc}}=\delta(p,a_{J}^{c})$.
Therefore $p^{'}\sim p$.

\noindent
(v)
If $i\ne j$, then $\Lambda(p_{i})\ne \Lambda(p_{j})$ by (iv).
Hence $\Lambda(p_{i})\cap  \Lambda(p_{j})=\emptyset$.
Therefore $\Lambda(p_{1})\sqcup\cdots \sqcup\Lambda(p_{M})\subset \Lambda$.
By (iii) and the decomposability of the branching function $f^{(\sigma)}$, 
$\Lambda(p_{1})\sqcup\cdots \sqcup\Lambda(p_{M})=\Lambda$.
This implies the statement.
\qedh

\noindent
{\it Proof of Theorem \ref{Thm:mainone}.}
Assume that $J=(j_{i})_{i=1}^{k}\in\nset{k}$. When $l=1$, $Q_{J}=\{q_{0}\}$. 
Let $J_{\sigma^{-1}}\equiv(\sigma^{-1}(j_{1}),\ldots, \sigma^{-1}(j_{k}))$.
Then we can check that $\lambda(q_{0},a_{J})=b_{J_{\sigma^{-1}}}$ and
$P(J)\circ \psi_{\sigma}=P(J_{\sigma^{-1}})$ independently.
Hence the assertion is verified.
Assume that $l\geq 2$.
By applying the correspondence between branching function systems and
permutative representations,
we see that the decomposition in Lemma \ref{lem:direct} (v)
implies that in (\ref{eqn:simpdeco}).
By definition of $J_{i}$ and applying Lemma \ref{lem:direct} (i), (ii)
to each component in the decomposition, 
the statement holds.
\qedh

\noindent
By Theorem \ref{Thm:mainone}, 
it is not necessary for computation of branching law (\ref{eqn:simpdeco})
to prepare any representation space. 
Further Theorem \ref{Thm:mainone} implies the following:
%
%
\begin{prop}
\label{prop:connected}
If the Mealy diagram of $\mms$ has $M$ connected components,
then $P(J)\circ \psi_{\sigma}$ has $M$ components 
of direct sum at least for each $J$.
\end{prop}

%
%
\sftt{Examples}
\label{section:fifth}
We show examples of permutative endomorphism of $\con$ and compute 
their branching laws by using the Mealy diagram according to 
Theorem \ref{Thm:mainone}.
Recall $E_{N,l}$ in (\ref{eqn:edef}).
Here we often denote $(j_{1},\ldots,j_{k})$ by $j_{1}\cdots j_{k}$ simply.
%
%
\ssft{$E_{2,2}$}
\label{subsection:etwo}
In \cite{PE01}, we show that there are $16$ equivalence classes in $E_{2,2}$
and there are $5$ irreducible and proper classes ${\cal E}$ in them.
We treat $3$ elements in ${\cal E}$ here. 
For each $\sigma\in {\goth S}_{2,2}$,
$\mms=(Q,\Sigma,\Delta,\delta,\lambda)$ consists of $Q=\{q_{1},q_{2}\}$,
$\Sigma=\{a_{1},a_{2}\}$ and $\Delta=\{b_{1},b_{2}\}$.

Define a transposition $\sigma\in {\goth S}_{2,2}$ by $\sigma(1,1)\equiv (1,2)$.
Then $\psi_{\sigma}$ and the Mealy diagram ${\cal D}(\mms)$ of 
$\mms$ are as follows:

\noindent
%
\def\bcircle#1{\circle{300}\put(-50,-30){$q_{#1}$}}
\def\labee#1#2{\put(0,0){$a_{#1}/b_{#2}$}}
\def\numbers{\put(-800,0){\bcircle{1}}\put(800,0){\bcircle{2}}}
\def\edgethree{
\thicklines
\qbezier(-650,-30)(0,-200)(650,-30)
\qbezier(-650,30)(0,200)(650,30)
\put(-650,-20){\vector(-3,1){0}}
\put(670,20){\vector(3,-1){0}}
}
\def\edges{\thicklines\put(0,0){\edgethree}}
\def\semi{
\qbezier(150,0)(0,150)(-150,0)
\qbezier(150,0)(200,-75)(150,-150)
\qbezier(-150,0)(-200,-75)(-150,-150)
\put(-110,-170){\vector(1,-1){0}}}
\def\labels{
\put(-1470,-40){\labee{2}{1}}
\put(1160,-40){\labee{2}{2}}
\put(-170,-250){\labee{1}{2}}
\put(-170,150){\labee{1}{1}}
}
%
\setlength{\unitlength}{.024713mm}
\begin{picture}(3011,620)(-3099,-300)
\thicklines
\put(0,0){\labels}
\put(0,0){\numbers}
\put(0,0){\edges}
%
\put(850,10){\rotatebox{-90}{\semi}}
\put(-1075,-20){\rotatebox{90}{\semi}}
\put(-3500,-50){$\left\{
\begin{array}{ll}
\psi_{\sigma}(s_{1})\equiv s_{1}s_{2}s_{1}^{*}+s_{1}s_{1}s_{2}^{*},\\
&\\
\psi_{\sigma}(s_{2})\equiv s_{2},\\
\end{array}
\right.$}
\end{picture}

\noindent
$\psi_{\sigma}$ is irreducible and proper (Table II in \cite{PE01}).
We denote $\psi_{\sigma}$ by $\psi_{12}$ in convenience.
We show several branching laws by $\psi_{12}$:
{\footnotesize
\[
\begin{array}{c|c|c|c}
\mbox{input}&\mbox{cycles}&\mbox{outputs}
&\mbox{branching law}\\
\hline
a_{1}&q_{1}q_{2}&b_{1}b_{2}&P(1)\circ 
\psi_{12}=
P(12)\\
\hline
a_{2}&q_{1},q_{2}&b_{1},b_{2}
&P(2)\circ \psi_{12}=
P(1)\oplus P(2)\\
\hline
a_{1}a_{2}&q_{1}q_{2}q_{2}q_{1}
&b_{1}b_{2}b_{2}b_{1}
&P(12)\circ \psi_{12}=
P(1122)\\
\hline
a_{1}a_{1}a_{2}a_{2}
&q_{1}q_{2}q_{1}q_{1},q_{2}q_{1}q_{2}q_{2}
&b_{1}b_{2}b_{1}b_{1},b_{2}b_{1}b_{2}b_{2}
&P(1122)\circ \psi_{12}=P(1112)\oplus P(1222)\\
\end{array}
\]
}

\noindent
Focusing attention on closed paths in ${\cal D}(\textsf{M}_{\sigma})$, 
we can verify the following:
%
%
\begin{prop}
\label{prop:psonetwo}
For each $J\in\{1,2\}^{*}_{1}$,
there are $J_{1},J_{2}$ or $J_{3}$ such that
\[
P(J)\circ \psi_{12}=
\left\{
\begin{array}{ll}
P(J_{1})\oplus P(J_{2})\quad &
(n_{1}(J)=\mbox{ even} ),\\
&\\
P(J_{3})\quad &(n_{1}(J)= \mbox{ odd})\\
\end{array}
\right.
\]
where $n_{1}(J)\equiv \sum_{l=1}^{k}(2-j_{l})$
for $J=(j_{1},\ldots,j_{k})\in\{1,2\}^{k}$.
\end{prop}

Let $\sigma\in {\goth S}_{2,2}$ be
a transposition defined by $\sigma(1,1)\equiv (2,1)$.
Then $\psi_{\sigma}$, ${\cal D}(\mms)$
and branching laws of $\psi_{\sigma}$ are given as follows:

\noindent
%
\def\bcircle#1{\circle{300}\put(-50,-30){$q_{#1}$}}
\def\labee#1#2{\put(0,0){$a_{#1}/b_{#2}$}}
\def\numbers{\put(-800,0){\bcircle{1}}\put(800,0){\bcircle{2}}}
\def\edgethree{
\thicklines
\qbezier(-650,-30)(0,-200)(650,-30)
\qbezier(-650,30)(0,200)(650,30)
\put(-650,-20){\vector(-3,1){0}}
\put(670,20){\vector(3,-1){0}}
}
\def\edges{\thicklines\put(0,0){\edgethree}}
\def\semi{
\qbezier(150,0)(0,150)(-150,0)
\qbezier(150,0)(200,-75)(150,-150)
\qbezier(-150,0)(-200,-75)(-150,-150)
\put(-110,-170){\vector(1,-1){0}}}
\def\labels{
\put(-1470,-40){\labee{1}{2}}
\put(1160,-40){\labee{2}{2}}
\put(-170,-250){\labee{1}{1}}
\put(-170,150){\labee{2}{1}}
}
%
\setlength{\unitlength}{.024713mm}
\begin{picture}(3011,500)(-3099,-200)
\thicklines
\put(0,0){\labels}
\put(0,0){\numbers}
\put(0,0){\edges}
%
\put(850,10){\rotatebox{-90}{\semi}}
\put(-1075,-20){\rotatebox{90}{\semi}}
\put(-3500,-50){$\left\{
\begin{array}{ll}
\psi_{\sigma}(s_{1})\equiv s_{2}s_{1}s_{1}^{*}+s_{1}s_{2}s_{2}^{*},\\
&\\
\psi_{\sigma}(s_{2})\equiv s_{1}s_{1}s_{1}^{*}+s_{2}s_{2}s_{2}^{*},\\
\end{array}
\right.$}
\end{picture}

{\small
\[
\begin{array}{c|c|c|c}
\mbox{input}&\mbox{cycles}&\mbox{outputs}
&\mbox{branching law}\\
\hline
a_{1}&q_{1}&b_{2}&P(1)\circ\psi_{\sigma}=P(2)\\
\hline
a_{2}&q_{2}&b_{2}&P(2)\circ\psi_{\sigma}=P(2)\\
\hline
a_{1}a_{2}&q_{2}q_{1}&b_{1}b_{2}&P(12)\circ\psi_{\sigma}=P(11)\\
\hline
a_{1}a_{1}a_{2}&
q_{2}q_{1}q_{1}&b_{1}b_{2}b_{1}&P(112)\circ\psi_{\sigma}=P(112)\\
\hline
a_{1}a_{2}a_{2}&
q_{2}q_{1}q_{2}&b_{1}b_{1}b_{2}&P(122)\circ\psi_{\sigma}=P(112)\\
\end{array}
\]
}

Let $\sigma\in {\goth S}_{2,2}$ be defined by
$\sigma(1,1)\equiv (2,2)$,
$\sigma(1,2)\equiv (1,1)$,
$\sigma(2,1)\equiv (2,1)$,
$\sigma(2,2)\equiv (1,2)$.
Then $\psi_{\sigma}$, ${\cal D}(\mms)$
and branching laws are as follows:

\noindent
%
\def\bcircle#1{\circle{300}\put(-50,-30){$q_{#1}$}}
\def\labee#1#2{\put(0,0){$a_{#1}/b_{#2}$}}
\def\numbers{\put(-800,0){\bcircle{1}}\put(800,0){\bcircle{2}}}
\def\edgethree{
\thicklines
\qbezier(-650,-30)(0,-200)(650,-30)
\qbezier(-650,30)(0,200)(650,30)
\put(-650,-20){\vector(-3,1){0}}
\put(670,20){\vector(3,-1){0}}
\qbezier(-750,-130)(0,-500)(750,-130)
\qbezier(-750,130)(0,500)(750,130)
\put(-740,-110){\vector(-1,1){0}}
\put(740,120){\vector(1,-1){0}}
}
\def\edges{\thicklines\put(0,0){\edgethree}}
\def\semi{
\qbezier(150,0)(0,150)(-150,0)
\qbezier(150,0)(200,-75)(150,-150)
\qbezier(-150,0)(-200,-75)(-150,-150)
\put(-110,-170){\vector(1,-1){0}}}
\def\labels{
\put(-170,-250){\labee{1}{2}}
\put(-170,150){\labee{2}{2}}
\put(-170,-420){\labee{2}{1}}
\put(-170,350){\labee{1}{1}}
}
%
\setlength{\unitlength}{.024713mm}
\begin{picture}(3011,920)(-3499,-290)
\thicklines
\put(0,0){\labels}
\put(0,0){\numbers}
\put(0,0){\edges}
\put(-3500,0){$\left\{
\begin{array}{ll}
\psi_{\sigma}(s_{1})\equiv s_{2}s_{2}s_{1}^{*}+s_{1}s_{1}s_{2}^{*},\\
&\\
\psi_{\sigma}(s_{2})\equiv s_{2}s_{1}s_{1}^{*}+s_{1}s_{2}s_{2}^{*},\\
\end{array}
\right.$}
\end{picture}

{\small
\[
\begin{array}{c|c|c|c}
\mbox{input}&\mbox{cycles}&\mbox{outputs}
&\mbox{branching law}\\
\hline
a_{1}&q_{1}q_{2}&b_{1}b_{2}&P(1)\circ 
\psi_{\sigma}=
P(12)\\
\hline
a_{2}&q_{1}q_{2}&b_{2}b_{1}&P(2)\circ 
\psi_{\sigma}=
P(12)\\
\hline
a_{1}a_{2}&q_{1}q_{2},q_{2}q_{1}
&b_{1}b_{1},b_{2}b_{2}&P(12)\circ 
\psi_{\sigma}=
P(11)\oplus  P(22)\\
\end{array}
\]
}
%
%
\ssft{$E_{3,2}$}
\label{subsection:fifthtwo}
Note that $\#E_{2,2}=2^{2}!=24$ and $\#E_{3,2}=3^{2}!\sim 3.6\times 10^{5}$.
Hence it is difficult to classify every element in 
$E_{3,2}$ by computing its branching laws in comparison with the case $E_{2,2}$.
We see that $\mms=(\{q_{1},q_{2},q_{3}\},
\{a_{1},a_{2},a_{3}\}$, $\{b_{1},b_{2},b_{3}\},
\delta,\lambda)$ for each $\sigma\in {\goth S}_{3,2}$.
$\rho_{\nu}$ in (\ref{eqn:nakeq}) belongs to $E_{3,2}$.

Let $\sigma\in {\goth S}_{3,2}$ be a transposition by $\sigma(1,1)\equiv (1,2)$. 
Then $\psi_{\sigma}$, ${\cal D}(\mms)$ 
and branching laws are as follows:

\noindent
%
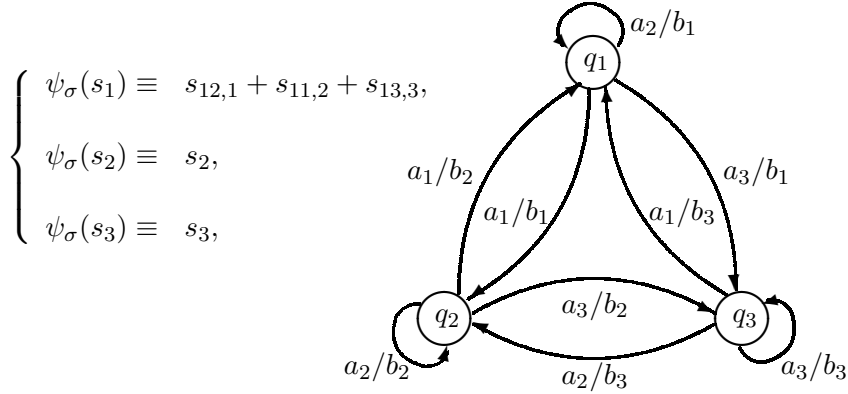
\def\bcircle#1{\circle{300}\put(-50,-20){$q_{#1}$}}
\def\labee#1#2{\put(0,0){$a_{#1}/b_{#2}$}}
\def\numbers{
\put(-800,0){\bcircle{2}}
\put(0,1380){\bcircle{1}}
\put(800,0){\bcircle{3}}
}
\def\edgethree{
\thicklines
\qbezier(-650,-30)(0,-400)(650,-30)
\qbezier(-650,30)(0,400)(650,30)
\put(-650,-20){\vector(-2,1){0}}
\put(650,30){\vector(2,-1){0}}
}
\def\edges{
\thicklines
\put(0,0){\edgethree}
\put(230,700){\rotatebox{120}{\edgethree}}
\put(-560,690){\rotatebox{-120}{\edgethree}}
}
\def\semi{
\qbezier(150,0)(0,150)(-150,0)
\qbezier(150,0)(200,-75)(150,-150)
\qbezier(-150,0)(-200,-75)(-150,-150)
\put(-110,-170){\vector(1,-1){0}}}
\def\labels{
\put(200,1550){\labee{2}{1}}
\put(-1340,-300){\labee{2}{2}}
\put(1010,-300){\labee{3}{3}}
%
\put(-1000,750){\labee{1}{2}}
\put(700,750){\labee{3}{1}}
\put(-170,-350){\labee{2}{3}}
\put(-170,30){\labee{3}{2}}
\put(-590,530){\labee{1}{1}}
\put(300,530){\labee{1}{3}}
}
%
\setlength{\unitlength}{.024713mm}
\begin{picture}(3011,1960)(-2899,-300)
\thicklines
\put(0,0){\labels}
\put(0,0){\numbers}
\put(0,0){\edges}
%
\put(-20,1620){\semi}
\put(850,-100){\rotatebox{-120}{\semi}}
\put(-1055,-130){\rotatebox{120}{\semi}}
\put(-3155,830){$
\left\{
\begin{array}{rl}
\psi_{\sigma}(s_{1})\equiv& s_{12,1}
+s_{11,2}+s_{13,3},\\
&\\
\psi_{\sigma}(s_{2})\equiv &s_{2},\\
&\\
\psi_{\sigma}(s_{3})\equiv &s_{3},\\
\end{array}
\right.
$}
\end{picture}

{\small
\[
\begin{array}{c|c|c|c}
\mbox{input}&\mbox{cycles}&\mbox{outputs}
&\mbox{branching law}\\
\hline
a_{1}&q_{1}q_{2}&b_{1}b_{2}&
P(1)\circ \psi_{\sigma}=
P(12)\\
\hline
a_{2}&q_{1},q_{2}&b_{1},b_{2}
&P(2)\circ \psi_{\sigma}=
P(1)\oplus P(2)\\
\hline
a_{3}&q_{3}&b_{3}
&P(3)\circ \psi_{\sigma}=
P(3)\\
\end{array}
\]
}

\noindent
where $s_{ij,k}\equiv s_{i}s_{j}s_{k}^{*}$.
From this, we see that $\psi_{\sigma}^{n}$
is proper and irreducible for each $n\geq 1$,
and $\psi_{\sigma}$ and $\rho_{\nu}$ are not equivalent.
%
%
\ssft{$E_{4,2}$}
\label{subsection:fourththree}
Define $\sigma\in {\goth S}_{4,2}$ by

\noindent
{\small
\[
\begin{array}{c|cccccccccccccccc}
J&11&12&13&14&21&22&23&24&31&32&33&34&41&42&43&44\\
\hline
\sigma(J)&11&21&31&41&12&22&43&42&32&23&13&33&44&24&14&34\\
\end{array}
\]
}

\noindent
Then $\psi_{\sigma}$ and ${\cal D}(\mms)$ are as follows:
\[
\begin{array}{ll}
\psi_{\sigma}(s_{1})\equiv 
s_{11,1}+s_{21,2}+s_{31,3}+s_{41,4},
&
\psi_{\sigma}(s_{2})\equiv 
s_{12,1}+s_{22,2}+s_{43,3}+s_{42,4},\\
&\\
\psi_{\sigma}(s_{3})\equiv 
s_{32,1}+s_{23,2}+s_{13,3}+s_{33,4},
&
\psi_{\sigma}(s_{4})\equiv 
s_{44,1}+s_{24,2}+s_{14,3}+s_{34,4},\\
\end{array}
\]

\noindent
%
%
%
\def\semi{
\qbezier(150,0)(0,150)(-150,0)
\qbezier(150,0)(200,-75)(150,-150)
\qbezier(-150,0)(-200,-75)(-150,-150)
\put(-110,-170){\vector(1,-1){0}}
}
\def\bcircle#1{\circle{400}\put(-70,-30){$q_{#1}$}}
\def\labee#1#2{\put(-70,0){$a_{#1}/b_{#2}$}}
\def\numbertwo{
\put(2400,500){\bcircle{2}}
\put(2390,820){\semi}
\put(1990,480){\rotatebox{90}{\semi}}
\put(2430,180){\rotatebox{180}{\semi}}
\put(2550,520){\rotatebox{270}{\semi}}
\put(2250,980){\labee{1}{1}}
\put(1610,450){\labee{2}{2}}
\put(2890,450){\labee{4}{4}}
\put(2290,-80){\labee{3}{3}}
}
\def\numberone{
\put(-1390,-100){\labee{1}{1}}
\put(-170,-250){\labee{2}{2}}
\put(-600,0){\bcircle{1}}
\put(-1010,-20){\rotatebox{90}{\semi}}
\put(-500,-270){\rotatebox{210}{\semi}}
}
\def\numberfour{
\put(1090,-100){\labee{1}{1}}
\put(-80,-250){\labee{2}{2}}
\put(600,0){\bcircle{4}}
\put(760,20){\rotatebox{-90}{\semi}}
\put(410,-290){\rotatebox{-210}{\semi}}
}
\def\numberthree{
\put(260,1280){\labee{1}{1}}
\put(-20,1320){\semi}
\put(0,1000){\bcircle{3}}
}
\def\sideedge{
\qbezier(-660,1500)(-2500,1100)(-1380,50)
\qbezier(-660,1400)(-1500,800)(-1280,50)
\qbezier(-550,1260)(-500,200)(-1150,0)
}
\def\leftedge{
\put(0,0){\sideedge}
\put(-640,1410){\vector(2,1){0}}
\put(-640,1510){\vector(3,1){0}}
\put(-1160,-20){\vector(-2,-1){0}}
}
\def\rightedge{
\put(10,10){\reflectbox{\sideedge}}
\put(550,1280){\vector(0,1){0}}
\put(1270,50){\vector(-1,-2){0}}
\put(1380,50){\vector(-1,-1){0}}
}
\def\figure{
\put(-610,500){\labee{2}{3}}
\put(-1700,1100){\labee{3}{3}}
\put(-1190,780){\labee{4}{4}}
\put(1330,1100){\labee{3}{3}}
\put(890,780){\labee{4}{4}}
\put(-150,10){\labee{4}{4}}
\put(330,500){\labee{3}{2}}
}
%
%
\setlength{\unitlength}{.020227mm}
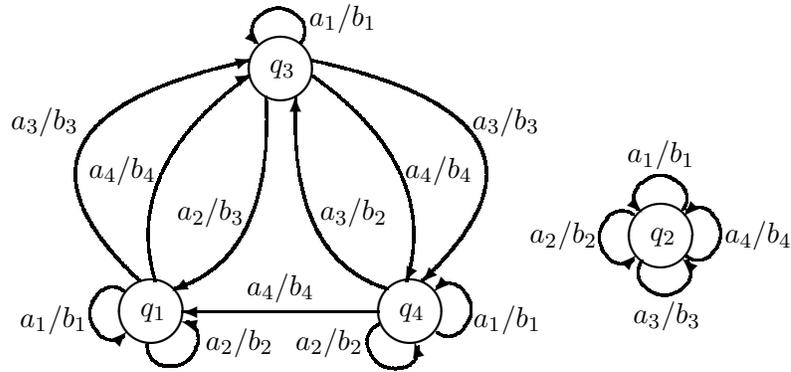
\begin{picture}(2756,2406)(-1799,-500)
\thicklines
\put(0,0){\figure}
\put(-250,-100){\numberone}
\put(100,-100){\numbertwo}
\put(0,500){\numberthree}
\put(250,-100){\numberfour}
\put(650,-100){\vector(-1,0){1300}}
\put(450,50){\leftedge}
\put(-450,40){\rightedge}
\end{picture}

\noindent
When $J=(1)$, $\delta(q_{i},a_{1})=q_{i}$
and $\lambda(q_{i},a_{1})=b_{1}$ for each $i=1,2,3,4$.
Therefore $P(1)\circ \psi_{\sigma} =P(1)\oplus P(1)\oplus P(1)\oplus P(1)$.
In the same way, we have
\[P(2)\circ \psi_{\sigma}=P(2)\oplus P(2)\oplus P(2), \quad
P(4)\circ \psi_{\sigma}=P(4)\oplus P(444).\]
This is an example of Proposition \ref{prop:connected}.

%
%
\ssft{Canonical endomorphism}
\label{subsection:fourthfour}
The Mealy diagram associated with the canonical endomorphism
$\rho$ of $\con$ (see $\S$\ref{section:third}) is given as follows:

\noindent
%
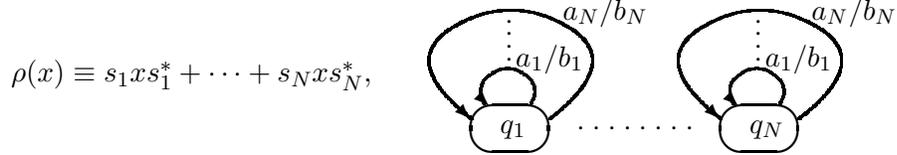
\def\bcircle#1{\oval(500,300)\put(-60,-30){$q_{#1}$}}
\def\labee#1#2{\put(0,0){$a_{#1}/b_{#2}$}}
\def\numbers{\put(-800,0){\bcircle{1}}\put(800,0){\bcircle{N}}}
\def\semi{
\qbezier(250,300)(0,400)(-250,300)
\qbezier(250,300)(800,75)(250,-350)
\qbezier(-250,300)(-800,75)(-250,-350)
\put(-250,-350){\vector(1,-1){0}}}
\def\semitwo{
\qbezier(150,0)(0,150)(-150,0)
\qbezier(150,0)(200,-75)(150,-150)
\qbezier(-150,0)(-200,-75)(-150,-150)
\put(-110,-170){\vector(1,-1){0}}
}
\def\labels{\multiput(-370,-40)(100,0){8}{$\cdot$}}
\def\crowns{
\put(800,310){\semi}
\put(780,210){\semitwo}
\multiput(770,300)(0,80){4}{$\cdot$}
\put(840,300){\labee{1}{1}}
\put(1140,600){\labee{N}{N}}
}
%
\setlength{\unitlength}{.020713mm}
\begin{picture}(3011,1000)(-3899,-150)
\thicklines
\put(0,0){\labels}
\put(0,0){\numbers}
%
\put(0,100){\crowns}
\put(-1600,100){\crowns}
\put(-4000,300){$\rho(x)\equiv
s_{1}xs_{1}^{*}+\cdots
+s_{N}xs_{N}^{*}$,}
\end{picture}

\noindent
In this case, there is no transition among different states.
We see that $P(J)\circ\rho=P(J)^{\oplus N}$
for each $J\in\nset{*}_{1}$ where
$P(J)^{\oplus N}$ is the direct sum of $N$ copies of $P(J)$.
In general, $\pi\circ \rho=\pi^{\oplus N}$
for any representation $\pi$ of $\con$.

%
%
\ssft{$E_{2,3}$}
\label{subsection:fourthfive}
Let $\sigma\in {\goth S}_{2,3}$ be 
a transposition by $\sigma(1,1,1)\equiv (1,2,1)$.
Then $\psi_{\sigma}\in E_{2,3}$, ${\cal D}(\mms)$ 
and branching laws are as follows:

\[\left\{
\begin{array}{l}
\psi_{\sigma}(s_{1})\equiv 
s_{121}s_{11}^{*}+s_{112}s_{12}^{*}+s_{111}s_{21}^{*}+s_{122}s_{22}^{*},\\
\\
\psi_{\sigma}(s_{2})\equiv s_{2},
\end{array}
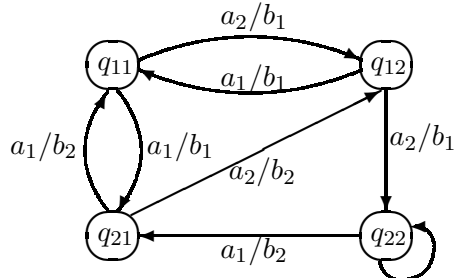
\right.\]

\noindent
%
\def\bcircle#1{\oval(300,280)\put(-90,-30){$q_{#1}$}}
\def\labee#1#2{\put(0,0){$a_{#1}/b_{#2}$}}
\def\numbers{
\put(-800,0){\bcircle{21}}
\put(800,0){\bcircle{22}}
\put(-800,1000){\bcircle{11}}
\put(800,1000){\bcircle{12}}
}
\def\edgethree{
\thicklines
\qbezier(-650,-30)(0,-300)(650,-30)
\qbezier(-650,30)(0,300)(650,30)
\put(-650,-20){\vector(-2,1){0}}
\put(650,30){\vector(2,-1){0}}
}
\def\edgefour{
\thicklines
\qbezier(-350,-30)(0,-300)(350,-30)
\qbezier(-350,30)(0,300)(350,30)
\put(-350,-30){\vector(-2,1){0}}
\put(350,40){\vector(2,-1){0}}
}
\def\edges{
\thicklines
\put(0,1000){\edgethree}
\put(-970,490){\rotatebox{90}{\edgefour}}
\put(650,0){\vector(-1,0){1300}}
\put(800,860){\vector(0,-1){720}}
\put(-690,120){\vector(2,1){1460}}
}
\def\semi{
\qbezier(150,0)(0,150)(-150,0)
\qbezier(150,0)(200,-75)(150,-150)
\qbezier(-150,0)(-200,-75)(-150,-150)
\put(-110,-170){\vector(1,-1){0}}}
\def\labels{
\put(-1400,480){\labee{1}{2}}
\put(-600,480){\labee{1}{1}}
\put(700,-390){\labee{2}{2}}
%
\put(-170,880){\labee{1}{1}}
\put(-170,1240){\labee{2}{1}}
\put(-170,-130){\labee{1}{2}}
\put(-120,330){\labee{2}{2}}
\put(820,530){\labee{2}{1}}
}
%
\setlength{\unitlength}{.022713mm}
\begin{picture}(3011,1420)(-2000,-100)
\thicklines
\put(0,0){\labels}
\put(0,0){\numbers}
\put(0,0){\edges}
%
\put(860,-150){\rotatebox{-135}{\semi}}
\end{picture}

{\small
\[
\begin{array}{c|c|c|c}
\mbox{input}&\mbox{cycles}&\mbox{outputs}
&\mbox{branching law}\\
\hline
a_{1}&q_{11}q_{21}
&b_{1}b_{2}&
P(1)\circ \psi_{\sigma}=P(12)\\
\hline
a_{2}&q_{22}
&b_{2}&
P(2)\circ \psi_{\sigma}=P(2)\\
\hline
a_{1}a_{2}&q_{12}q_{11}
&b_{1}b_{1}&
P(12)\circ \psi_{\sigma}=P(11)\\
\hline
a_{1}a_{1}a_{2}
&q_{12}q_{11}q_{21}
&b_{1}b_{1}b_{2}
&P(112)\circ \psi_{\sigma}=P(112)\\
\end{array}
\]
}
We see that $\psi_{\sigma}^{n}$ is irreducible and proper for each $n\geq 1$.

Let $\sigma\in {\goth S}_{2,3}$ be defined by
the product $\sigma=\sigma^{'}\circ \sigma^{''}$ of two transpositions
$\sigma^{'}$ and $\sigma^{''}$ defined by $\sigma^{'}(1,1,1)\equiv (1,2,1)$ 
and $\sigma^{''}(1,1,2)\equiv (1,2,2)$, respectively.
In this case $\psi_{\sigma}=\psi_{12}\in E_{2,2}$ in $\S$\ref{subsection:etwo}.
${\cal D}(\mms)$ is as follows:

\noindent
%
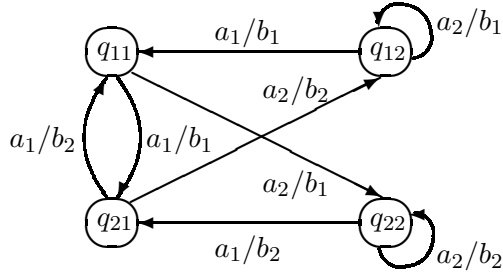
\def\bcircle#1{\oval(300,280)\put(-90,-20){$q_{#1}$}}
\def\labee#1#2{$a_{#1}/b_{#2}$}
\def\numbers{
\put(-800,0){\bcircle{21}}
\put(800,0){\bcircle{22}}
\put(-800,1000){\bcircle{11}}
\put(800,1000){\bcircle{12}}
}
\def\edgefour{
\thicklines
\qbezier(-350,-30)(0,-300)(350,-30)
\qbezier(-350,30)(0,300)(350,30)
\put(-350,-30){\vector(-2,1){0}}
\put(350,40){\vector(2,-1){0}}
}
\def\edges{
\thicklines
\put(-970,490){\rotatebox{90}{\edgefour}}
\put(650,0){\vector(-1,0){1300}}
\put(650,1000){\vector(-1,0){1300}}
\put(-680,880){\vector(2,-1){1460}}
\put(-690,120){\vector(2,1){1460}}
}
\def\semi{
\qbezier(150,0)(0,150)(-150,0)
\qbezier(150,0)(200,-75)(150,-150)
\qbezier(-150,0)(-200,-75)(-150,-150)
\put(-110,-170){\vector(1,-1){0}}}
\def\labels{
\put(-1400,440){\labee{1}{2}}
\put(-600,440){\labee{1}{1}}
\put(1100,-240){\labee{2}{2}}
%
\put(-200,1070){\labee{1}{1}}
%
\put(-200,-200){\labee{1}{2}}
\put(80,730){\labee{2}{2}}
\put(80,180){\labee{2}{1}}
%
\put(1100,1130){\labee{2}{1}}
}
%
\setlength{\unitlength}{.022713mm}
\begin{picture}(3011,1720)(-2099,-300)
\thicklines
\put(0,0){\labels}
\put(0,0){\numbers}
\put(0,0){\edges}
%
\put(860,-150){\rotatebox{-135}{\semi}}
\put(840,1180){\rotatebox{-45}{\semi}}
\end{picture}

\noindent
We can verify that branching laws of $\psi_{\sigma}$
coincide with those of $\psi_{12}$.
%
%
\ssft{$E_{2,4}$}
\label{subsection:fourthsix}
Define a transposition $\sigma\in {\goth S}_{2,4}$ 
by $\sigma(1,1,1,1)\equiv (1,2,1,1)$.
Then $\psi_{\sigma}\in E_{2,4}$, ${\cal D}(\mms)$
and branching laws are given as follows:
{\small
\[
\!\psi_{\sigma}(s_{1})\equiv 
s_{1211}s_{111}^{*}+
s_{1112}s_{112}^{*}+
s_{112}s_{12}^{*}+
s_{1111}s_{211}^{*}+
s_{1212}s_{212}^{*}+
s_{122}s_{22}^{*},\quad
\psi_{\sigma}(s_{2})\equiv s_{2},
\]
}

\noindent
%
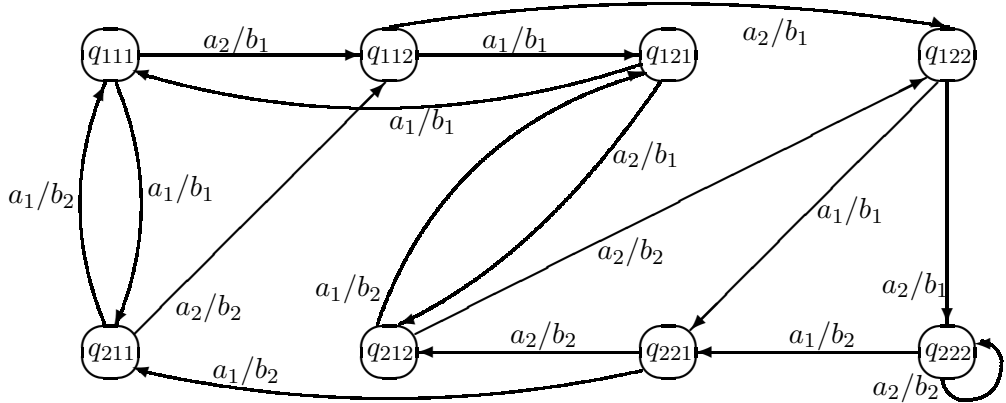
\def\bcircle#1{\oval(300,280)\put(-120,-20){$q_{#1}$}}
\def\labee#1#2{$a_{#1}/b_{#2}$}
\def\unin{1500}
\def\unitwo{3000}
\def\unithree{4500}
\def\numbers{
\put(0,1600){\bcircle{111}}
\put(\unin,1600){\bcircle{112}}
\put(\unitwo,1600){\bcircle{121}}
\put(\unithree,1600){\bcircle{122}}
\put(0,0){\bcircle{211}}
\put(\unin,0){\bcircle{212}}
\put(\unitwo,0){\bcircle{221}}
\put(\unithree,0){\bcircle{222}}
}
\def\edgethree{
\thicklines
\qbezier(-650,-30)(0,-300)(650,-30)
\qbezier(-650,30)(0,300)(650,30)
\put(-650,-20){\vector(-2,1){0}}
\put(650,30){\vector(2,-1){0}}
}
\def\edges{
\thicklines
\put(-170,790){\rotatebox{90}{\edgethree}}
\put(150,1600){\vector(1,0){1200}}
\put(1650,1600){\vector(1,0){1200}}
\qbezier(1500,1750)(3000,2000)(4500,1750)
\put(4500,1750){\vector(3,-1){0}}
\qbezier(150,1500)(1500,1100)(2850,1550)
\put(120,1520){\vector(-3,1){0}}
\qbezier(1550,150)(2300,500)(2950,1450)
\put(1550,150){\vector(-3,-1){0}}
\put(4450,1450){\vector(-1,-1){1330}}
\put(4500,1450){\vector(0,-1){1300}}
\put(130,100){\vector(1,1){1350}}
\qbezier(1430,150)(1800,1200)(2870,1500)
\put(2870,1510){\vector(3,1){0}}
\put(1640,100){\vector(2,1){2750}}
\put(2850,0){\vector(-1,0){1200}}
\qbezier(150,-100)(1500,-400)(2850,-100)
\put(120,-80){\vector(-3,1){0}}
\put(4350,0){\vector(-1,0){1200}}
}
\def\semi{
\qbezier(150,0)(0,150)(-150,0)
\qbezier(150,0)(200,-75)(150,-150)
\qbezier(-150,0)(-200,-75)(-150,-150)
\put(-110,-170){\vector(1,-1){0}}}
\def\labels{
\put(200,840){\labee{1}{1}}
\put(500,1640){\labee{2}{1}}
\put(2000,1640){\labee{1}{1}}
\put(3400,1680){\labee{2}{1}}
\put(1500,1200){\labee{1}{1}}
\put(2700,1020){\labee{2}{1}}
\put(3800,720){\labee{1}{1}}
\put(4150,320){\labee{2}{1}}
\put(-550,800){\labee{1}{2}}
\put(350,200){\labee{2}{2}}
\put(1100,300){\labee{1}{2}}
\put(2620,500){\labee{2}{2}}
\put(550,-140){\labee{1}{2}}
\put(2150,50){\labee{2}{2}}
\put(3650,50){\labee{1}{2}}
\put(4100,-230){\labee{2}{2}}
}
%
\setlength{\unitlength}{.024713mm}
\begin{picture}(3011,2320)(100,-200)
\thicklines
\put(0,0){\labels}
\put(0,0){\numbers}
\put(0,0){\edges}
%
\put(4570,-150){\rotatebox{-135}{\semi}}
\end{picture}

{\small
\[
\begin{array}{c|c|c|c}
\mbox{input}&\mbox{cycles}&\mbox{outputs}
&\mbox{branching law}\\
\hline
a_{1}&q_{111}q_{211}
&b_{1}b_{2}&
P(1)\circ \psi_{\sigma}=P(12)\\
\hline
a_{2}&q_{222}
&b_{2}&
P(2)\circ \psi_{\sigma}=P(2)\\
\hline
a_{1}a_{2}&q_{212}q_{121}
&b_{2}b_{1}&
P(12)\circ \psi_{\sigma}=P(12)\\
\hline
a_{1}a_{1}a_{2}&q_{112}q_{121}q_{111}
&b_{1}b_{1}b_{1}&
P(112)\circ \psi_{\sigma}=P(111)\\
\end{array}
\]
}

\noindent
{\bf Acknowledgement:}
The author would like to thank Takeshi Nozawa for useful comment on this article.


\end{document}